\documentclass[reqno,12pt]{amsart}
\usepackage{amsmath,amssymb,color}
\usepackage{amsthm}
\usepackage{comment}
\usepackage{graphicx}

\newtheorem{theorem}{Theorem}
\newtheorem{corollary}{Corollary}

\newtheorem{remark}{Remark}

\def\re{\mathbb{R}}

\def\eps{\varepsilon}

\def\pd{\partial}
\def\ol{\overline}

\def\la{\lambda}

\def\({\left(}
\def\){\right)}

\def\pd{\partial}

\def\|{\Vert}

\begin{document}
\title[On a nonlocal overdetermined problem]{A bifurcation analysis on a nonlocal overdetermined problem}

\author{Kazuki Sato}

\address{
Department of Mathematics, Osaka Metropolitan University \\
3-3-138, Sumiyoshi-ku, Sugimoto-cho, Osaka, Japan \\
}

\email{sf22817a@st.omu.ac.jp \\}

\author{Futoshi Takahashi}

\address{
Department of Mathematics, Osaka Metropolitan University \\
3-3-138, Sumiyoshi-ku, Sugimoto-cho, Osaka, Japan \\
}

\email{futoshi@omu.ac.jp \\}

\begin{abstract}
In this paper, we study an overdetermined problem with Kirchhoff type nonlocal terms related to the celebrated work by Serrin.
We obtain the precise number of solutions according to the value of the bifurcation parameter 
and study asymptotics of bifurcation curves of solutions when the bifurcation parameter is large in some cases.
\end{abstract}


\date{\today}

\dedicatory{}

\maketitle

\section{Introduction}

Let $\Omega \subset \re^N$ be a bounded domain with $C^2$ boundary.
In this paper, we consider the following overdetermined elliptic problem with Kirchhoff type nonlocal terms:
\begin{equation}
\label{OD_AB}
	\begin{cases}
	&-A\(\| u \|_{L^{p_1}(\Omega)}, \| \nabla u \|_{L^{q_1}(\Omega)} \) \Delta u = \la B\( \| u \|_{L^{p_2}(\Omega)}, \| \nabla u \|_{L^{q_2}(\Omega)}\) \quad \text{in} \ \Omega, \\
	&u = 0 \quad \text{on} \ \pd\Omega, \\
	&\frac{\pd u}{\pd \nu} = c \quad \text{on} \ \pd\Omega,
	\end{cases}
\end{equation}
where $p_1, p_2, q_1, q_2 > 0$, $\la > 0$ are given constants, 
$c$ is an unknown constant, $\nu$ is an outer unit normal to $\pd\Omega$,
and $A(s, t)$, $B(s, t)$ are positive continuous functions in $(s, t) \in \re_+ \times \re_+$.
The positive constant $\la$ plays a role as a bifurcation parameter.

In a celebrated paper by Serrin \cite{Serrin}, see also Weinberger \cite{Weinberger}, 
it is proven that if the overdetermined elliptic problem
\begin{equation}
\label{OD}
	\begin{cases}
	&-\Delta u = 1 \quad \text{in} \ \Omega, \\
	&u = 0 \quad \text{on} \ \pd\Omega, \\
	&\frac{\pd u}{\pd \nu} = c \quad \text{on} \ \pd\Omega
	\end{cases}
\end{equation}
admits a solution $u \in C^2(\ol{\Omega})$, then $\Omega$ must be a ball and $u$ is radially symmetric. 
Thus if $\Omega = B_R(x_0)$, a ball of radius $R$ with center $x_0$, then the unique solution of \eqref{OD} must be of the form
\begin{equation}
\label{U_R}
	U_{R,x_0}(x) = \frac{R^2-|x-x_0|^2}{2N}.
\end{equation}
By a simple computation, we have 
\begin{align*}
	&\| U_{R, x_0} \|^p_{L^p(B_R(x_0))} = \(\frac{1}{2N}\)^p \(\frac{\omega_{N-1}}{2}\) R^{2p+N} B\(\frac{N}{2}, p+1\), \\ 
	&\| \nabla U_{R, x_0} \|^q_{L^q(B_R(x_0))} = \(\frac{1}{N}\)^q \(\frac{\omega_{N-1}}{N+q}\) R^{N+q} 
\end{align*}
for any $p, q > 0$, where $\omega_{N-1}$ is the area of the unit sphere in $\re^N$ and $B(x, y) = \int_0^1 t^{x-1} (1-t)^{y-1} dt$ denotes the Beta function. 

After \cite{Serrin} \cite{Weinberger}, many researches have been done around this topic and various generalizations (differential operators, spaces, and so on) have been obtained up to now.
We refer to a survey paper \cite{NT} and the references there in.
However, to the best of authors' knowledge, overdetermined problems with Kirchhoff type nonlocal terms have not been studied so far.
Thus we initiate to study it in this paper.

\vspace{1em}
In this paper, following an argument in \cite{ACM}, we prove
\begin{theorem}
\label{Theorem:P}
Let $A, B \in C(\re_+ \times \re_+ ; \re_+)$.
Then if \eqref{OD_AB} admits a solution $u \in C^2(\ol{\Omega})$, then $\Omega$ must be a ball and $u$ must be radially symmetric. 

On the other hand, consider the system of equations 
with respect to $(s_1, s_2, t_1, t_2) \in \re_+ \times \re_+ \times \re_+ \times \re_+$:
\begin{equation}
\label{System}
	\begin{cases}
	&s_1 = \la \dfrac{B(s_2, t_2)}{A(s_1, t_1)} \| U_{R, x_0} \|_{L^{p_1}(B_R(x_0))} \\
	&s_2 = \la \dfrac{B(s_2, t_2)}{A(s_1, t_1)} \| U_{R, x_0} \|_{L^{p_2}(B_R(x_0))} \\
	&t_1 = \la \dfrac{B(s_2, t_2)}{A(s_1, t_1)} \| \nabla U_{R, x_0} \|_{L^{q_1}(B_R(x_0))} \\
	&t_2 = \la \dfrac{B(s_2, t_2)}{A(s_1, t_1)} \| \nabla U_{R, x_0} \|_{L^{q_2}(B_R(x_0))}
	\end{cases}
\end{equation}
where $U_{R, x_0}$ defined in \eqref{U_R} is the unique solution of \eqref{OD} for $\Omega = B_R(x_0)$.
Then for $\la > 0$, the problem \eqref{OD_AB} for $\Omega = B_R(x_0)$ has the same number of solutions of the system of equations \eqref{System}.
Also the number of solutions of \eqref{System} is the same as the number of solutions to the equation
\begin{equation}
\label{Single}
	g(s) = \la \| U_{R, x_0} \|_{L^{p_1}(B_R(x_0))}
\end{equation}
with respect to $s > 0$, where
\begin{equation}
\label{g(s)}
	g(s) =  \frac{s A\(s, \dfrac{\| \nabla U_{R, x_0} \|_{L^{q_1}(B_R(x_0))}}{\| U_{R, x_0} \|_{L^{p_1}(B_R(x_0))}} s \)}{B\(\dfrac{\| U_{R, x_0} \|_{L^{p_2}(B_R(x_0))}}{\| U_{R, x_0} \|_{L^{p_1}(B_R(x_0))}}s, \dfrac{\| \nabla U_{R, x_0} \|_{L^{q_2}(B_R(x_0))}}{\| U_{R, x_0} \|_{L^{p_1}(B_R(x_0))}}s \)}
\end{equation}
 
Moreover, any solution $u_{\la}$ of \eqref{OD_AB} is of the form 
\begin{equation}
\label{u_form}
	u_{\la}(x) = s_1 \frac{U_{R,x_0}(x)}{\| U_{R, x_0} \|_{L^{p_1}(B_R(x_0))}}
\end{equation}
where $s_1$ is a solution of \eqref{Single}.
\end{theorem}

\vspace{1em}
As a result of Theorem \ref{Theorem:P}, many bifurcation diagrams can be obtained for specific $A$ and $B$ in \eqref{OD_AB}.
Next are typical examples of the result.


\begin{corollary}
\label{Cor2}
Let $\Omega = B_R(x_0)$ and let $A(s, t) = (s-a)^2 + b$, $B(s,t) = s$ in \eqref{OD_AB} for $a, b > 0$.
Then 
\begin{enumerate}
\item[(i)] if $0 < \la < \dfrac{b}{\| U_{R,x_0} \|_{L^{p_2}(B_R(x_0))}}$, there is no solution to \eqref{OD_AB}, 
\item[(ii)] if $\la = \dfrac{b}{\| U_{R,x_0} \|_{L^{p_2}(B_R(x_0))}}$, there is a unique solution to \eqref{OD_AB}, 
\item[(iii)] if $\dfrac{b}{\| U_{R,x_0} \|_{L^{p_2}(B_R(x_0))}} < \la < \dfrac{(a^2 +b)}{\| U_{R,x_0} \|_{L^{p_2}(B_R(x_0))}}$, there are just two solutions to \eqref{OD_AB}, 
\item[(iv)] if $\la \ge \dfrac{(a^2 + b)}{\| U_{R,x_0} \|_{L^{p_2}(B_R(x_0))}}$, there is a unique solution to \eqref{OD_AB}.
\end{enumerate}
Also any solution $u_{\la}$ is of the form
\[
	u_{\la}(x) = s_{\la} \frac{U_{R,x_0}(x)}{\| U_{R,x_0} \|_{L^{p_1}(B_R(x_0))}}, \quad x \in B_R(x_0)
\]
where $s_{\la}$ is a positive solution of the quadratic equation
\[
	(s-a)^2 + b = \la \| U_{R,x_0} \|_{L^{p_2}(B_R(x_0))}.
\]
\end{corollary}
%
%

\begin{corollary}
\label{Cor4}
Let $\Omega = B_R$ and let $A(s, t) = e^s$, $B(s, t) = t^r$ in \eqref{OD_AB} where $r > 1$.
Then 
\begin{enumerate}
\item[(i)] 
if $0 < \la \| U_{R,x_0} \|_{L^{p_1}(B_R(x_0))} < \(\frac{e}{r-1}\)^{r-1} \(\frac{\| U_{R,x_0} \|_{L^{p_1}(B_R(x_0))}}{\| \nabla U_{R,x_0} \|_{L^{q_2}(B_R(x_0))}}\)^r$, there is no solution to \eqref{OD_AB}, 
\item[(ii)] 
if $\la \| U_{R,x_0} \|_{L^{p_1}(B_R(x_0))} = \(\frac{e}{r-1}\)^{r-1} \(\frac{\| U_{R,x_0} \|_{L^{p_1}(B_R(x_0))}}{\| \nabla U_{R,x_0} \|_{L^{q_2}(B_R(x_0))}}\)^r$, there is a unique solution to \eqref{OD_AB}, 
\item[(iii)] 
if $\la \| U_{R,x_0} \|_{L^{p_1}(B_R(x_0))} > \(\frac{e}{r-1}\)^{r-1} \(\frac{\| U_{R,x_0} \|_{L^{p_1}(B_R(x_0))}}{\| \nabla U_{R,x_0} \|_{L^{q_2}(B_R(x_0))}}\)^r$, there are just two solutions $u_{1,\la}$, $u_{2,\la}$ to \eqref{OD_AB}.
\end{enumerate}
Also we have
\begin{align*}
	&u_{1,\la}(x) = \frac{\la^{-\frac{1}{r-1}}}{\| \nabla U_{R,x_0} \|_{L^{q_2}(B_R(x_0))}^{\frac{r}{r-1}}} \\
	&\times \( 1 + \frac{1}{r-1}(1 + o(1)) \la^{-\frac{1}{r-1}} \frac{\| U_{R,x_0} \|_{L^{p_1}(B_R(x_0))}}{\| \nabla U_{R,x_0} \|_{L^{q_2}(B_R(x_0))}^{\frac{r}{r-1}}} \) U_{R,x_0}(x) 
\end{align*}
for $x \in B_R(x_0)$ as $\la \to \infty$,
and
\[
	u_{2,\la}(x) = \left\{ \log \la + (r-1) (\log \log \la)(1 + o(1)) \right\} \frac{U_{R,x_0}(x)}{\| U_{R,x_0} \|_{L^{p_1}(B_R(x_0))}} 
\]
for $x \in B_R(x_0)$ as $\la \to \infty$.
\end{corollary}

\vspace{1em}
Here we mention some related works.
Let $I = (-1, 1) \subset \re$ be an interval.
T. Shibata performed a bifurcation analysis on the one-dimensional Dirichlet problems with Kirchhoff type nonlocal terms of the form
\[
	\begin{cases}
	&-A\(\| u \|_{L^{q_1}(I)}, \| u' \|_{L^{r_1}(I)} \) u''(x) = \la B\( \| u \|_{L^{q_2}(I)}, \| u' \|_{L^{r_2}(I)}\) u^p(x) \quad x \in I,\\
	&u(x) > 0 \quad x \in I, \\
	&u(\pm 1) = 0
	\end{cases}
\]
where $p > 1$ and $\la > 0$.
We refer the readers to a series of works by T. Shibata \cite{Shibata(JMAA)}, \cite{Shibata(BVP)}, \cite{Shibata(ANONA)}, \cite{Shibata(QTDS)}. 
Also, in the same spirit, we performed a bifurcation analysis on the nonlocal boundary blowup problem in one space dimension
\[
	\begin{cases}
	&A\(\| u \|_{L^{q_1}(I)}, \| u' \|_{L^{r_1}(I)} \) u''(x) = \la B\( \| u \|_{L^{q_2}(I)}, \| u' \|_{L^{r_2}(I)}\) u^p(x) \quad x \in I,\\
	&u(x) > 0 \quad x \in I, \\
	&\lim_{x \to \pm 1} u(x) = +\infty
	\end{cases}
\]
where $p > 1$, $\la > 0$, $q_1, q_2 \in (0, \frac{p-1}{2})$, and $r_1, r_2 \in (0, \frac{p-1}{p+1})$,
see \cite{Inaba-TF}, \cite{Sato-TF}.
We consider the problem \eqref{OD_AB} is a higher dimensional analogue of the problems above.

\section{Proof of Theorem \ref{Theorem:P}.}

In this section, we prove Theorem \ref{Theorem:P}.

\begin{proof}
First assume that there exists a solution $u \in C^2(\ol{\Omega})$ of \eqref{OD_AB} and put 
\[
	v = \gamma u
\] 
where $\gamma > 0$ is chosen so that
\[
	\gamma^{-1} = \la \frac{B(\| u \|_{L^{p_2}(\Omega)}, \| \nabla u \|_{L^{q_2}(\Omega)})}{A(\| u \|_{L^{p_1}(\Omega)}, \| \nabla u \|_{L^{q_1}(\Omega)})} 
\]
Then
\begin{align*}
	-\Delta v(x) = \gamma (-\Delta u(x)) \overset{\eqref{OD}}{=} \gamma \la \dfrac{B(\| u \|_{L^{p_2}(\Omega)}, \| \nabla u \|_{L^{q_2}(\Omega)})}{A(\| u \|_{L^{p_1}(\Omega)}, \| \nabla u \|_{L^{q_1}(\Omega)})} = 1
\end{align*}
in $\Omega$.
Also we see $v = 0$ on $\pd\Omega$ and $\frac{\pd v}{\pd \nu} = \gamma c$ is a constant on $\pd\Omega$.
Thus by the result of Serrin, we see that $\Omega$ must be a ball, say $\Omega = B_R(x_0)$ for some $R > 0$ and $x _0 \in \re^N$ and 
$v \equiv U_{R,x_0}(x)$. 
This implies that 
\begin{equation}
\label{u_form}
	u(x) = \gamma^{-1} U_{R,x_0}(x), \quad \nabla u(x) = \gamma^{-1} \nabla U_{R,x_0}(x).
\end{equation}
Define
\[
	s_1 = \| u \|_{L^{p_1}(B_R(x_0))}, \ s_2 = \| u \|_{L^{p_2}(B_R(x_0))}, \ t_1 = \| \nabla u \|_{L^{q_1}(B_R(x_0))}, \ t_2 = \| \nabla u \|_{L^{q_2}(B_R(x_0))}.
\]
Then by \eqref{u_form}, we have
\[
	\begin{cases}
	&s_1 = \gamma^{-1} \| U_{R,x_0} \|_{L^{p_1}(B_R(x_0))}, \\
	&s_2 = \gamma^{-1} \| U_{R,x_0} \|_{L^{p_2}(B_R(x_0))}, \\
	&t_1 = \gamma^{-1} \| \nabla U_{R,x_0} \|_{L^{q_1}(B_R(x_0))}, \\
	&t_2 = \gamma^{-1} \| \nabla U_{R,x_0} \|_{L^{q_2}(B_R(x_0))},
	\end{cases}
\]
which is equivalent to \eqref{System}.
This shows that 
\[
	(s_1, s_2, t_1, t_2) = (\| u \|_{L^{p_1}(B_R(x_0))}, \| u \|_{L^{p_2}(B_R(x_0))}, \| \nabla u \|_{L^{q_1}(B_R(x_0))}, \| \nabla u \|_{L^{q_2}(B_R(x_0))})
\]
is a solution to \eqref{System} and thus
\[
	\sharp \{ u : \text{solutions of \eqref{OD_AB}} \} \le \sharp \{ (s_1, s_2, t_1, t_2) \in (\re_+)^4 : \text{solutions of \eqref{System}} \},
\]
where $\sharp A$ denotes the cardinality of the set $A$.

On the other hand, let $\Omega = B_R(x_0)$ for some $R > 0$ and $x_0 \in \re^N$ 
and let $(s_1, s_2, t_1, t_2) \in (\re_+)^4$ be any solution to \eqref{System}.
Note that by \eqref{System}, we see
\begin{align}
\label{System2}
	&\frac{s_1}{\| U_{R, x_0} \|_{L^{p_1}(B_R(x_0))}} = \frac{s_2}{\| U_{R, x_0} \|_{L^{p_2}(B_R(x_0))}} = \frac{t_1}{\| \nabla U_{R,x_0} \|_{L^{q_1}(B_R(x_0))}} = \frac{t_2}{\| \nabla U_{R, x_0} \|_{L^{q_2}(B_R(x_0))}} \notag \\
	&= \la \frac{B(s_2, t_2)}{A(s_1, t_1)}.
\end{align}
Thus if we define 
\begin{align*}
	u(x) &= s_1 \frac{U_{R, x_0}(x)}{\| U_{R, x_0} \|_{L^{p_1}(B_R(x_0))}} \\
	&\(= s_2 \frac{U_{R, x_0}(x)}{\| U_{R, x_0} \|_{L^{p_2}(B_R(x_0))}} = t_1 \frac{U_{R, x_0}(x)}{\| \nabla U_{R, x_0} \|_{L^{q_1}(B_R(x_0))}} = t_2 \frac{U_{R, x_0}(x)}{\| \nabla U_{R, x_0} \|_{L^{q_2}(B_R(x_0))}} \),
\end{align*}
then we have $u(x) > 0$, $u = 0$, $\frac{\pd u}{\pd \nu} = const.$ on $\pd B_R(x_0)$ and 
\[
	s_1 = \| u \|_{L^{p_1}(B_R(x_0))}, \ s_2 = \| u \|_{L^{p_2}(B_R(x_0))}, \ t_1 = \| \nabla u \|_{L^{q_1}(B_R(x_0))}, \ t_2 = \| \nabla u \|_{L^{q_2}(B_R(x_0))}.
\]
Moreover, by the definition of $u$, we have
\begin{align*}
	-A\(\| u \|_{L^{p_1}(B_R(x_0))}, \| \nabla u \|_{L^{q_1}(B_R(x_0))}\) & \Delta u(x) = -A(s_1, t_1) \Delta u(x) \\
	&= -A(s_1, t_1) \frac{s_1}{\| U_{R, x_0} \|_{L^{p_1}(B_R(x_0))}} \underbrace{\Delta U_{R, x_0}(x)}_{=-1} \\
	&= A(s_1, t_1) \frac{s_1}{\| U_{R, x_0} \|_{L^{p_1}(B_R(x_0))}} \\
	&\overset{\eqref{System2}}= \la B(s_2, t_2).
\end{align*}
This shows that
\[
	\sharp \{ u : \text{solutions of \eqref{OD_AB}} \} \ge \sharp \{ (s_1, s_2, t_1, t_2) \in (\re_+)^4 : \text{solutions of \eqref{System}} \}.
\]
Thus the number of solutions of \eqref{OD_AB} for $\Omega = B_R(x_0)$ and that of \eqref{System} are the same.

Also by \eqref{System2}, we can rewrite the system of equations \eqref{System} into a single equation for $s = s_1$ 
\[
	s_1 =
	\la \dfrac{B\(\dfrac{\| U_{R, x_0} \|_{L^{p_2}(B_R(x_0))}}{\| U_{R, x_0} \|_{L^{p_1}(B_R(x_0))}} s_1, \dfrac{\| \nabla U_{R, x_0} \|_{L^{q_2}(B_R(x_0))}}{\| U_{R, x_0} \|_{L^{p_1}(B_R(x_0))}} s_1 \)}
{A\(s_1, \dfrac{\| \nabla U_{R, x_0} \|_{L^{q_1}(B_R(x_0))}}{\| U_{R, x_0} \|_{L^{p_1}(B_R(x_0))}} s_1 \)} \| U_{R, x_0} \|_{L^{p_1}(B_R(x_0))}
\]
which is equivalent to \eqref{Single} with $g(s)$ in \eqref{g(s)}.
Thus the number of solutions of \eqref{System} for $\Omega = B_R(x_0)$ and that of \eqref{Single} are also the same. 
\end{proof}

\begin{remark}
Theorem \ref{Theorem:P} can be generalized to the overdetermined problem of the form
\[
	\begin{cases}
	&-A\(\| u \|_{\vec{p}}, \| \nabla u \|_{\vec{q}} \) \Delta u = \la B\( \| u \|_{\vec{r}}, \| \nabla u \|_{\vec{s}} \), \quad \text{in} \ \Omega, \\
	&u = 0, \quad \text{on} \ \pd\Omega, \\
	&\frac{\pd u}{\pd \nu} = c, \quad \text{on} \ \pd\Omega
	\end{cases}
\]
where
\begin{align*}
	&A\(\| u \|_{\vec{p}}, \| \nabla u \|_{\vec{q}} \) = A\(\| u \|_{L^{p_1}(\Omega)}, \cdots, \| u \|_{L^{p_I}(\Omega)}, \| \nabla u \|_{L^{q_1}(\Omega)}, \cdots, \| \nabla u \|_{L^{q_J}(\Omega)}\), \\
	&B\(\| u \|_{\vec{r}}, \| \nabla u \|_{\vec{s}} \) = B\(\| u \|_{L^{r_1}(\Omega)}, \cdots, \| u \|_{L^{r_K}(\Omega)}, \| \nabla u \|_{L^{s_1}(\Omega)}, \cdots, \| \nabla u \|_{L^{s_L}(\Omega)}\)
\end{align*}
for $\{ p_i \}_{i=1}^I, \{ q_j \}_{j=1}^J, \{ r_k \}_{k=1}^K, \{ s_l \}_{l=1}^L \subset \re_+$.
We treat $A, B$ of the forms in Theorem \ref{Theorem:P} just for a simple presentation.
\end{remark}

\section{Proof of corollaries.}

In this section, we prove corollaries in \S 1.
By Theorem \ref{Theorem:P}, the number of solutions of \eqref{OD_AB} is the same as the number of solutions of \eqref{Single} when $\Omega = B_R(x_0)$,
and we have the expression of solution $u(x) = s \frac{U_{R,x_0}(x)}{\| U_{R,x_0} \|_{L^{p_1}}}$ to \eqref{OD_AB} where $s > 0$ is any solution of \eqref{Single}. 
Here and in what follows, $\| \cdot \|_p$ stands for $\| \cdot \|_{L^p(B_R(x_0))}$ for any $p > 0$.
Thus we just need to solve the equation \eqref{Single} with $\la > 0$ for specific $A$ and $B$.


\vspace{1em}
{\it Proof of Corollary \ref{Cor2}}:\hspace{1em}
Since $A(s, t) = \{ (s-a)^2 + b \}$ and $B(s, t) = s$, the function $g(s)$ defined in \eqref{g(s)} becomes
\begin{align*}
	g(s) = \frac{\| U_{R, x_0} \|_{p_1}}{\| U_{R, x_0} \|_{p_2}} \left\{ (s - a)^2 + b \right\}.
\end{align*}
Thus $g(s)$ is a quadratic function with respect to $s > 0$ and 
\[
	\min_{s > 0} g(s) = g(a) = b \frac{\| U_{R, x_0} \|_{p_1}}{\| U_{R, x_0} \|_{p_2}}, \quad g(0) = (a^2+b) \frac{\| U_{R, x_0} \|_{p_1}}{\| U_{R, x_0} \|_{p_2}}.
\]
Now, the equation \eqref{Single} reads 
\[
	\frac{\| U_{R, x_0} \|_{p_1}}{\| U_{R, x_0} \|_{p_2}} \left\{ (s - a)^2 + b \right\} = \la \| U_{R, x_0} \|_{p_1}.
\]
By a consideration of the graph of $g(s)$, we see that the number of positive solutions $s$ of \eqref{Single} according to the value of $\la > 0$ is:
\begin{enumerate}
\item[(i)] $0$, if $0 < \la \| U_{R, x_0} \|_{p_1} < g(a) = \min_{s > 0} g(s)$,
\item[(ii)] $1$, if $\la \| U_{R, x_0} \|_{p_1} = g(a)$,
\item[(iii)] $2$, if $g(a) < \la \| U_{R, x_0} \|_{p_1} < g(0)$,
\item[(iv)] $1$, if $\la \| U_{R, x_0} \|_{p_1} \ge g(0)$.
\end{enumerate}
This completes the proof of Corollary \ref{Cor2}.
\qed


\vspace{1em}
{\it Proof of Corollary \ref{Cor4}}:\hspace{1em}
Let $A(s, t) = e^s$ and $B(s, t) = t^r$ for $r > 1$. 
Then the function $g(s)$ defined in \eqref{g(s)} becomes
\[
	g(s) = \( \frac{\| U_{R, x_0} \|_{p_1}}{\| \nabla U_{R, x_0} \|_{q_2}} \)^r e^s s^{1-r}.
\]
A computation yields that 
\[
	\min_{s > 0} g(s) = g(r-1) = \( \frac{\| U_{R, x_0} \|_{p_1}}{\| \nabla U_{R, x_0} \|_{q_2}} \)^r \(\frac{e}{r-1}\)^{r-1} 
\]
and $\lim_{s \to +0} g(s) = \lim_{s \to \infty} g(s) = +\infty$.
Thus the number of solutions to \eqref{Single} is
\begin{enumerate}
\item[(i)] $0$, if $0 < \la \| U_{R, x_0} \|_{p_1} < \( \frac{\| U_{R, x_0} \|_{p_1}}{\| \nabla U_{R, x_0} \|_{q_2}} \)^r \(\frac{e}{r-1}\)^{r-1}$,
\item[(ii)] $1$, if $\la \| U_{R, x_0} \|_{p_1} = \( \frac{\| U_{R, x_0} \|_{p_1}}{\| \nabla U_{R, x_0} \|_{q_2}} \)^r \(\frac{e}{r-1}\)^{r-1}$,
\item[(iii)] $2$, if $\la \| U_{R, x_0} \|_{p_1} > \( \frac{\| U_{R, x_0} \|_{p_1}}{\| \nabla U_{R, x_0} \|_{q_2}} \)^r \(\frac{e}{r-1}\)^{r-1}$.
\end{enumerate}
Let $0 < s_{1, \la} < s_{2, \la}$ be two solutions of $g(s) = \la \| U_{R, x_0} \|_{p_1}$ in the case (iii).
Then we see $s_{1,\la} \to 0$ and $s_{2,\la} \to \infty$ as $\la \to \infty$.
Since
\[
	e^{s_{i,\la}} = s_{i,\la}^{r-1} \la K
\]
for $i=1,2$, where
\[
	K = \frac{\| \nabla U_{R, x_0} \|^r_{q_2}}{\| U_{R, x_0} \|^{r-1}_{p_1}}, 
\]
by taking $\log$ of the above equation, we obtain
\begin{equation}
\label{sila}
	s_{i, \la} = (r-1) \log s_{i,\la} + \log (\la K).
\end{equation}
Since $s_{1,\la} = o(1)$ as $\la \to \infty$, we have
\[
	o(1) = (r-1) \log s_{1,\la} + \log (\la K)
\]
as $\la \to \infty$.
By solving this with respect to $s_{1,\la}$, we see
\begin{equation}
\label{s1la}
	s_{1,\la} = (\la K)^{-\frac{1}{r-1}} e^{o(1)} = (\la K)^{-\frac{1}{r-1}} (1 + \delta)
\end{equation}
where $\delta \to 0$ as $\la \to \infty$.
Inserting \eqref{s1la} into the both sides of \eqref{sila} and computing, we have
\[
	(\la K)^{-\frac{1}{r-1}}(1 + \delta) = (r-1) \log ( 1 + \delta) = (r-1) \delta(1 + o(1)).
\]
From this, we have
\begin{align*}
	(\la K)^{-\frac{1}{r-1}} = (r-1) \frac{\delta}{1 + \delta} (1 + o(1)) = (r-1) \delta (1 + o(1)).
\end{align*}
Thus
\begin{align*}
	\delta = \frac{1}{r-1} (1 + o(1)) (\la K)^{-\frac{1}{r-1}} \quad (\la \to \infty).
\end{align*}
Inserting this into \eqref{s1la}, we have 
\[
	s_{1,\la} = (\la K)^{-\frac{1}{r-1}} \( 1 + \frac{1}{r-1}(1 + o(1)) (\la K)^{-\frac{1}{r-1}} \)
\]
and recalling $u_{1,\la} = s_{1,\la} \frac{U_{R, x_0}(x)}{\| U_{R, x_0} \|_{p_1}}$,
we obtain the formula for $u_{1,\la}$:
\begin{align*}
	u_{1,\la}(x) &= (\la K)^{-\frac{1}{r-1}} \( 1 + \frac{1}{r-1}(1 + o(1)) (\la K)^{-\frac{1}{r-1}} \) \frac{U_{R, x_0}(x)}{\| U_{R, x_0} \|_{p_1}} \\
	&= \frac{\la^{-\frac{1}{r-1}}}{\| \nabla U_{R, x_0} \|_{q_2}^{\frac{r}{r-1}}} \( 1 + \frac{1}{r-1}(1 + o(1)) \la^{-\frac{1}{r-1}} \frac{\| U_{R, x_0} \|_{p_1}}{\| \nabla U_{R, x_0} \|_{q_2}^{\frac{r}{r-1}}} \) U_{R, x_0}(x) 
\end{align*}
for $x \in B_R(x_0)$.

Similarly, from \eqref{sila}, we have
\[
	s_{2, \la} = (r-1) \log s_{2,\la} + \log (\la K).
\]
Now, since $s_{2, \la} \to \infty$ as $\la \to \infty$, we have
\[
	s_{2, \la} \(1 - (r-1) \underbrace{\frac{\log s_{2,\la}}{s_{2,\la}}}_{=o(1)} \) = \log (\la K) = (\log \la) \(1 + \underbrace{\frac{\log K}{\log \la}}_{= o(1)}\)
\]
as $\la \to \infty$.
From this, we have
\begin{equation}
\label{s2la}
	s_{2,\la} = \frac{1}{1-o(1)} (\log \la)(1 + o(1)) = (1+ \eps) \log \la
\end{equation}
as $\la \to \infty$, where $\eps \to 0$ as $\la \to \infty$. 
Again, inserting this into \eqref{sila}, we have
\[
	(1 + \eps) \log \la = (r-1) \left\{ \log (1 + \eps) + \log \log \la \right\} + \log \la + \log K,
\]
which leads to
\begin{align*}
	\eps &= (r-1)\frac{\log \log \la}{\log \la} \left\{ 1 + \underbrace{\frac{\log (1 + \eps)}{\log \log \la} + \frac{1}{r-1} \frac{\log K}{\log \log \la}}_{=o(1)} \right\} \\
	&= (r-1)\frac{\log \log \la}{\log \la} (1 + o(1))
\end{align*}
as $\la \to \infty$. Then coming back to \eqref{s2la} and recalling that $u_{2,\la} = s_{2,\la} \frac{U_{R, x_0}(x)}{\| U_{R, x_0} \|_{p_1}}$,
we obtain the asymptotic formula for $u_{2,\la}$:
\begin{align*}
	&s_{2, \la} = (1 + \eps) \log \la = \( 1 + (r-1)\frac{\log \log \la}{\log \la} (1 + o(1)) \) \log \la, \\
	&u_{2,\la}(x) = \left\{ \log \la + (r-1) (\log \log \la)(1 + o(1)) \right\} \frac{U_{R, x_0}(x)}{\| U_{R, x_0} \|_{p_1}}
\end{align*}
as $\la \to \infty$.
\qed

\vskip 0.5cm

\noindent\textbf{Acknowledgement.} 
The second author (F.T.) was supported by JSPS Grant-in-Aid for Scientific Research (B), No. 23H01084, 
and was partly supported by Osaka Central University Advanced Mathematical Institute (MEXT Joint Usage/Research Center on Mathematics and Theoretical Physics).

\end{document}